\date{} % %
\title{The Ising magnetization exponent on ${\mathbb Z}^2$ is $1/15$}

\author{Federico Camia \and Christophe Garban \and Charles M. Newman}

\documentclass[11pt,preprint]{article}

\usepackage[vmargin={3cm, 3.5cm},hmargin=2.5cm]{geometry}

\usepackage[pdfpagelabels]{hyperref}

\usepackage{latexsym} 
\usepackage{amsmath}
\usepackage{amssymb}
\usepackage{amsthm}
\usepackage{mathrsfs}
\usepackage{amsfonts}
\usepackage{graphicx}
\usepackage{url}

\usepackage{color}
\usepackage{bbm,verbatim}

\usepackage{minitoc} %% POUR LES TABLES MATIERES DANS LES CHAPITRES.

\input labelfig.tex

\usepackage{multirow} %pour les tableaux
\usepackage{array}
\newcolumntype{M}[1]{>{\centering}m{#1}}

\numberwithin{equation}{section}
\numberwithin{figure}{section}
\newtheorem{theorem}{Theorem}
\numberwithin{theorem}{section}
\newtheorem{corollary}[theorem]{Corollary}
\newtheorem{lemma}[theorem]{Lemma}

\newtheorem{proposition}[theorem]{Proposition}

\theoremstyle{remark}\newtheorem{remark}[theorem]{Remark}
\def\eqref#1{(\ref{#1})}
\let\qqed=\qed

\def\nn{\nonumber}

\def\QED{\qqed\medskip}
\let\qed=\QED

\newcommand{\Z}{\mathbb{Z}}
\newcommand{\N}{\mathbb{N}}

\def\SLEkk#1/{$\mathrm{SLE}(#1)$}
\def\SLEr#1/{$\mathrm{SLE(\kappa;#1)}$}
\def\SLEkr#1;#2/{$\mathrm{SLE(#1;#2)}$}
\def\SLEk/{\SLEkk{\kappa}/}
\def\SLEtwo/{\SLEkk2/}
\def\SLE/{$\mathrm{SLE}$}
\def\SLEab/{\SLEkr 4; {a/\hco-1}, {b/\hco-1}/}
\def\Var#1{\mathrm{Var}\big[ #1\big]}

\def \eps {\epsilon}
\def \P {\mathbb{P}}
\def \E {{\mathbb E}}

\def\md{\mid}

\def \p {{\partial}}

\def\FK#1#2#3{{\def\md{\bigm| } \P_{#1}^{\,#2}  \bigl[  #3 \bigr]}}
\def\EFK#1#2#3{{\def\md{\bigm| } \E_{#1}^{\,#2}  \bigl[  #3 \bigr]}}

\def\fk{\mathrm{FK}}
\def\free{\mathrm{free}}

\def \p {{\partial}}

\def\1{1\hspace{-2.55 mm}{1}}

\def\noopsort#1{}

\usepackage{mathrsfs}

\def\ni{\noindent}
\def\bi{\begin{itemize}}
\def\ei{\end{itemize}}

\def\ghost{\mathrm{ghost}}
\def\cluster{\mathcal{C}}

\def\<#1{\langle #1\rangle}
\def\g{\mathbf{g}}

\def\GHS{MR0266507}

\def\Wu{Wu}
\def\WuB{WuB}

\def\CGNSL{CGNSL}

\def\HVising{HVising} %??

\def\Wer09{MR2523462}

\def\SimonPphi2{SimonPphi2}

\def\Chayes{MR1678308}
\def\ABF87{ABF87}
\def\Griffith{Gr67}

\def\RSWfk{arXiv:0912.4253}
\def\GrimmettFK{MR2243761}

\def\GPS{MR2736153}

\def\Onsager{MR0010315}

\def\KestenScaling{MR88k:60174}

\def\SmirnovWerner{MR1879816}

\def\BeffaraDim6{MR2078552}

\begin{document}
\maketitle

\begin{abstract}
We prove that for the Ising model defined on the plane $\Z^2$ at $\beta=\beta_c$, the average magnetization under an external magnetic field $h>0$ behaves exactly like 
\[
\<{\sigma_0}_{\beta_c, h} \asymp h^{\frac 1 {15} }\,.
\]
The proof, which is surprisingly simple compared to an analogous result for percolation (i.e. that $\theta(p)=(p-p_c)^{5/36+o(1)}$ on the triangular lattice \cite{\SmirnovWerner,\KestenScaling}) relies on the GHS inequality as well as the  RSW theorem for FK percolation from \cite{\RSWfk}.  The use of GHS to obtain
inequalities involving critical exponents is not new; in this paper we show how it can be combined with RSW to
obtain matching upper and lower bounds for the average magnetization.
\end{abstract}

%\tableofcontents

\section{Introduction}

The classical Ising model on a finite domain $\Lambda_L:= [-L,L]^2 \subset \Z^2$ with $+$ boundary condition and
with external field $h\geq 0$ is  a probability measure on $\{-1,1\}^{\Lambda_L}$, $\P^{\beta,h,+}_{L}$,
defined as follows.
For any spin configuration $\sigma\in \{-1,1\}^{\Lambda_L}$, let 
\begin{equation}\label{}
E_L(\sigma):= -\sum_{x\sim y} \sigma_x \sigma_y - \sum_{x\in \p \Lambda_L} \sigma_x
\end{equation}
be the interaction energy,  where the first sum is over nearest neighbor pairs in $\Lambda_L$ and the second
is over sites in $\partial \Lambda_L$, the boundary of $\Lambda_L$. Let also 
\begin{equation}\label{}
M_L(\sigma):= \sum_{x\in \Lambda_L} \sigma_x
\end{equation}
be the total magnetization in $\Lambda_L$. The probability measure $\P^{\beta,h,+}_{L}$ on $\{-1,1\}^{\Lambda_L}$
is defined by 
\begin{equation}\label{ising-dist-h}
\FK{L}{\beta,h,+}{\sigma}:= \frac 1 {Z_{L, \beta, h}} e^{ -\beta\, E_L(\sigma) + h\, M_L(\sigma)}\,,
\end{equation}
where the partition function $Z_{L,\beta, h}$ is simply defined as $\sum_\sigma e^{ -\beta\, E_L(\sigma) + h\, M_L(\sigma)}$. 

As is well known, the measures $\P^{\beta,h,+}_{L}$ have a unique {\it infinite volume limit} as $L\to \infty$, that we denote by $\P_{\Z^2}^{\beta, h, +}$ (in fact, when $h>0$  or $\beta \leq \beta_c$, adding the $+$ boundary
condition does not have any effect on the infinite volume limit,  so the $+$ in the notation can be dropped).   
Since the work by Onsager (\cite{\Onsager}), it is known that if $h=0$, then this system undergoes a {\bf phase transition} at \begin{equation}\label{}
\beta_c = \frac 1 2 \ln(1+\sqrt{2})\,.
\end{equation}
See for example \cite{\GrimmettFK} and references therein and see \cite{\HVising} for a recent self-contained proof. 

%\medskip
In this paper, we are interested in the behavior of the Ising model on $\Z^2$ near its critical point at $\beta=\beta_c$ but with a small external magnetic field $h>0$. Our main theorem is the following,  where $\<{\cdot}_{\beta, h}$
denotes expectation with respect to $\P_{\Z^2}^{\beta, h, +}$ (when $h=0$, we sometimes drop the second index).
\begin{theorem}\label{th.hZ2}
Consider the Ising model on $\Z^2$ at $\beta_c$ with a positive external magnetic field 
$h>0$, then \footnote{In this paper $f(a) \asymp g(a)$ as $a \searrow 0$ means that $f(a)/g(a)$ is bounded
away from zero and $\infty$ while $f(a) \sim g(a)$ means that $f(a)/g(a) \to 1$ as $a \searrow 0$.}
\[
\<{\sigma_0}_{\beta_c, h} \asymp h^{\frac 1 {15} }\, .
\]
\end{theorem}

Our proof of Theorem \ref{th.hZ2} relies on the GHS inequality as well as the recent  RSW theorem
from \cite{\RSWfk}. While the use of GHS to obtain inequalities involving critical exponents is not new (see, e.g.,
Chapter 14 of \cite{FFS92} and references therein), we combine it in a novel way with RSW for Ising-FK percolation
to obtain, to the best of our knowledge, the first complete proof that the Ising magnetization exponent is $1/15$.

%\section{ Introduction}
\bigskip

Let us now give more background on the subject. The Ising model on the square lattice with no external magnetic field has been solved exactly,
and its solution has yielded a number of critical exponents. The first such exponent was obtained by Onsager.
\begin{theorem}[Onsager \cite{\Onsager}]\label{onsager}
 If we denote by $\langle \cdot \rangle_{\beta}^{+}$ expectation with respect to the
infinite volume measure with inverse temperature $\beta$, zero external magnetic field, and $+$
boundary condition, then
\begin{equation*}
\langle \sigma_0 \rangle_{\beta}^{+} \asymp |\beta - \beta_c|^{1/8}\,,
 \qquad \text{as }\beta \searrow \beta_c \, .
\end{equation*}
\end{theorem}

A second exponent follows from a celebrated result of T. T. Wu, which is crucial for the proof of our main 
result, Theorem \ref{th.hZ2}.
\begin{theorem}[T. T. Wu, see \cite{\WuB, \Wu}]\label{th.Wu}
There exists an explicit constant $c>0$ such that as $n \to \infty $
\begin{equation}\label{}
\rho(n):= %\<{\sigma_{(0,0)} \sigma_{(n,n)}}
\<{\sigma_{(0,0)} \sigma_{(n,n)}}_{\beta_c} \sim c\, n^{-1/4}\,.
\end{equation}
\end{theorem}
See Section 1.2 in \cite{\CGNSL} for a discussion of the status of Wu's result.  
See also Section~\ref{WWu} below 
where we state a theorem which does not assume Wu's result 
(in particular, it turns out that by assuming much less, one still obtains the exponent
$\frac 1 {15}$ in a weaker sense --- see Remark \ref{r.SLE}).
%\medskip

These critical exponents are defined for zero external magnetic field.
If one introduces a nonzero external magnetic field, the average magnetization is also nonzero,
and its behavior as we let the external field go to zero defines (modulo its existence proved in
this paper) another critical exponent,
$\delta$, via $\<{\sigma_0}_{\beta_c, h} \asymp h^{\frac 1 {\delta} }$.

The value $\delta=15$
 is suggested by non-rigorous scaling theory and can be understood heuristically by
considering the continuum scaling limit of the Ising 
model with a vanishing  external magnetic field
(i.e., the near-critical limit in which a nonzero $h$ scales appropriately to
zero), as we now
briefly explain. In the continuum scaling limit, the lattice is scaled by a factor $1/L$, with $L \to \infty$,
and one focuses, for example, on the magnetization in the unit square.

At the critical temperature and with zero external magnetic field, one can show \cite{CGNSL} that the random
variable $m_L := L^{-15/8} M_L(\sigma) = L^{-15/8} \sum_{x \in \Lambda_L} \sigma_x$ has a unique limit
in distribution as $L \to \infty$. The scaling factor $L^{-15/8}$ insures that the second moment of $m_L$ is bounded away from zero and infinity as $L \to \infty$. See \cite{CGNSL} for more details.

%At the critical temperature and with zero external magnetic field, one can show \cite{CGNSL} that the random
%variable $m_L := (L^2)^{-15/8} M_L(\sigma) = L^{-15/16} \sum_{x \in \Lambda_L} \sigma_x$ has a unique limit
%in distribution as $L \to \infty$. The scaling factor $(L^2)^{-15/8}$ 
%(asymptotically proportional to the number of spins in $\Lambda_L$ to
%the power $-15/8$) insures that the second moment of $m_L$ is bounded away from zero and infinity as $L \to \infty$,
%with the choice of the exponent $15/8$ dictated by the value $1/4$ of the two-point 
%function exponent in Theorem~\ref{th.Wu}.
%(Indeed, the scaling factor must scale like
%$\sqrt{\sum_{x,y \in {\Lambda_L}} \<{\sigma_x \sigma_y}_{\beta_c, 0}}$ as $L \to \infty$,
%see \cite{MR2504956} or \cite{CGNSL} for more details.)

Adding an external magnetic field $h$ yields the Ising distribution \eqref{ising-dist-h}, which contains the term
$h M_L(\sigma)$. In order to obtain a meaningful continuum scaling limit in this case, it appears necessary to
let $h$ scale as the scaling factor $L^{-15/8}$ discussed in the previous paragraph (see \cite{MR2504956}
or \cite{Camia} for a discussion of this point, and 
\cite{CGNproperties} for a construction of the continuum near-critical scaling
limit with a vanishing external magnetic field).

For any $1 \leq l < L$, we can write the expectation under $\P^{\beta,h,+}_{L}$ of the rescaled
magnetization $m_l(\sigma) = l^{-15/8} \sum_{x \in \Lambda_l} \sigma_x$ in $\Lambda_l$ as 
\begin{equation} \nonumber
\frac{\sum_{\sigma} m_l(\sigma) \exp(-\beta_c E_L(\sigma) + h M_L(\sigma)) } {\sum_{\sigma} \exp(-\beta_c E_L(\sigma) + h M_L(\sigma)) } \, .
%\stackrel{L\gg1}{\sim} a^{15/8} a^{-2} \langle S_0
%\rangle_{\beta_c,H=a^{15/8}} \sim a^{-1/8} (a^{15/8})^{1/15} = 1 \, .
\end{equation}
When $L \gg l \to \infty$, assuming the scaling law $\<{\sigma_0}_{\beta_c, h} \asymp h^{\frac 1 {\delta} }$ holds,
the above ratio behaves like 
\begin{equation} \nonumber
l^{-15/8} l^2 \langle \sigma_0 \rangle_{\beta_c,h=l^{-15/8}} \asymp l^{1/8} (l^{-15/8})^{1/\delta} \, .
\end{equation}
%If the choices of scaling factors discussed in the previous paragraphs are correct, and
If a nontrivial continuum scaling limit exists, as proved in \cite{CGNproperties}, one expects this quantity
to have a finite nonzero limit as $l \to \infty$, which requires $\delta=15$. (We note that, using the non-rigorous
scaling laws for critical exponents and the fact that the heat capacity exponent $\alpha$ is zero for the
two-dimensional Ising model, one can write $\delta = \frac{2-1/8}{1/8}$, where $1/8$ is
the order parameter exponent of Theorem \ref{onsager}---see, e.g., Section 16.3 of \cite{Huang87}.

\medskip

Some comments on the interpretation of our main result, Theorem \ref{th.hZ2}, are in order.
Since the quantity $\<{\sigma_0}_{\beta_c, h}$ can be interpreted as the probability that the origin is connected
to the {\it ghost vertex} in the 
appropriate FK percolation model defined on $\Z^2 \cup \{ \ghost \}$ (see Section \ref{s.lower}), one can think of
Theorem \ref{th.hZ2}, as an analog of the following theorem by Smirnov and Werner.

\begin{theorem}[Smirnov and Werner \cite{\SmirnovWerner}]\label{}
The density function $\theta(p)$ for site percolation on the triangular grid has the following behavior for $p>p_c=1/2$:
\begin{equation*}
\theta(p) = (p-1/2)^{5/36+o(1)}\,,
\end{equation*}
as $p\to 1/2+$. 
\end{theorem}

Theorem \ref{th.hZ2} is also an analog of the celebrated result of Onsager on the Ising spontaneous magnetization,
Theorem \ref{onsager} above, except that our result concerns the near-critical regime in the $h$ direction rather
than in the $\beta$ direction. 

\medskip

Let us end this introduction by stating the  %{\bf GHS inequality} 
 Griffiths-Hurst-Sherman inequality 
from \cite{\GHS} which will be essential to our work. 

\begin{theorem}[GHS inequality \cite{\GHS}]\label{th.GHS}
Let $G=(V,E)$ be a finite graph. Consider a ferromagnetic Ising model on this graph (i.e., the interactions $J_e$ for $e=\{i,j\} \in E$ are non-negative)
and assume furthermore that the external field 
$\mathbf{h} = (h_v)_{v\in V}$ (which may vary from one vertex to another) is  non-negative. Under such general assumptions, one has for any vertices 
$i,j,k\in V$:
\[
\langle \sigma_i\sigma_j \sigma_k\rangle - 
\Bigl( 
\langle \sigma_i\rangle\, \langle \sigma_j \sigma_k\rangle
+\langle \sigma_j\rangle\, \langle \sigma_i \sigma_k\rangle
+\langle \sigma_k\rangle\, \langle \sigma_i \sigma_j\rangle
\Bigr)
+2 \langle \sigma_i\rangle \langle \sigma_j\rangle \langle \sigma_k\rangle
\leq 0 \,.
\]
\end{theorem}

This inequality has the following useful corollary.
\begin{corollary}\label{cor.GHS}
Let $G=(V,E)$ be a finite graph and let $K \subset V$ be a non-empty subset of the vertices. Let us consider a ferromagnetic Ising model on $G$ with the spins in $K$ prescribed to be $+$ and with a constant magnetic field $h\geq 0$ on $V\setminus K$. Then the partition function of this model, i.e.,
\[
Z_{\beta, h} := \sum_{\sigma \in \{-,+\}^{V\setminus K}} \exp{\left(-\beta \sum_{i \sim j \in V} J_{\{i,j\}} \sigma_i \sigma_j + h \sum_{i \in V \setminus K} \sigma_i \right)}\,,
\]
satisfies
\[
\p^3_h \log(Z_{\beta,h}) \leq 0\,. 
\]
\end{corollary}

\noindent{\bf Proof.} Since the partition function $Z_{\beta,{\bf h}}$ for a (non-constant) external field $\bf h$ has
$\p_{h_i} \p_{h_j} \p_{h_k} \log Z_{\beta,{\bf h}}$ given by the LHS of the displayed inequality of Theorem~\ref{th.GHS},
this is an immediate corollary.
\QED

In the next section we will use the GHS inequality to obtain Eq. \eqref{e.1and2}; we will then use that bound,
combined with the inequalities in Proposition \ref{pr.1and2} below, to obtain an optimal upper bound for the
average magnetization with a nonzero external magnetic field. 

As pointed out to us by Hugo Duminil-Copin, the inequality \eqref{e.1and2} was already used in \cite{FFS92} (see Eq. (14.230),
p. 345) to obtain a lower bound for the spontaneous magnetization, leading to an inequality involving critical exponents.
The use of the GHS inequality to obtain bounds like Eq. \eqref{e.1and2} was apparently first proposed by the third
author in an unpublished 1982 preprint that was later included as an appendix of \cite{Newman86}. 
%In the next section we will use the bound in a different way and will combine it with the inequalities in Proposition
%\ref{pr.1and2} below to obtain an optimal upper bound for the average magnetization with a nonzero external
%magnetic field.

%The use of the GHS inequality to obtain bounds like Eq. \eqref{e.1and2} below
%was apparently first proposed by the third author in an unpublished preprint that was later included as an
%appendix of \cite{Newman86}. We note that Eq. \eqref{e.1and2} is used in \cite{FFS92} (see Eq. (14.230), p. 345)
%to obtain a lower bound for the spontaneous magnetization, leading to an inequality involving critical exponents.
%In the next section we will use the bound in a different way and will combine it with the inequalities in Proposition
%\ref{pr.1and2} below to obtain an optimal upper bound for the average magnetization with a nonzero external
%magnetic field.

\section{Proof of the upper bound}

Recall that, for any $L\in \N$, 
\[
M_L= \sum_{x\in \Lambda_L} \sigma_x
\]
is the (non-renormalized) magnetization in the square $\Lambda_L =[-L,L]^2$. If some boundary condition 
 $\eta$ is prescribed on
$\Lambda_L$, we will denote the magnetization by $M^\eta_L$. In the same fashion, we will denote by $M^{\Z^2}_L$ the magnetization in $\Lambda_L$
with boundary condition inherited from the full-plane $\Z^2$.  In particular, by translation invariance, one has 

\[
\<{\sigma_0}_{\beta_c,h}  = \frac 1 {|\Lambda_L|} %\EFK{}{\beta_c,h}{M^{\Z^2}_L}\,,
\<{M^{\Z^2}_L}_{\beta_c,h} \, ,
\]
for any $L\in \N$. 

By monotonicity one has, for any $L\geq 1$,
\begin{align}\label{e.dombymagn}
\<{\sigma_0}_{\beta_c,h}  \le \frac 1 {|\Lambda_L|} %\EFK{}{\beta_c,h,+}{M_L}\,.
\<{M_L^+}_{\beta_c,h} \, .
\end{align}
 
 For notational convenience, in the rest of the proof and in the next proposition,
we will use $\<{M_L}_{\beta_c,h,+}$ to denote $\<{M_L^+}_{\beta_c,h}$. The main idea
in the proof of the upper bound is to rewrite the expected magnetization
 $\<{M_L}_{\beta_c,h,+}$ %$\EFK{}{\beta_c,h,+}{M_L}$
as follows: 
\[
%\EFK{}{\beta_c,h,+}{M_L}
\<{M_L}_{\beta_c,h,+} = \frac 
{\<{M_L\, e^{h M_L}}_{\beta_c,0,+}}
{\<{e^{h M_L}}_{\beta_c,0,+}}
= \frac 
{\frac \p {\p h} \<{e^{h M_L}}_{\beta_c,0,+}}
{\<{e^{h M_L}}_{\beta_c,0,+}}\,,
\]
and then to apply the GHS inequality. Indeed the latter (recall Corollary \ref{cor.GHS}) says that,
for $+$ boundary conditions, 
\begin{eqnarray*}
& \frac{\p^3 }{\p h^3} \log \Bigl( \sum_\sigma e^{-\beta_c E_L(\sigma) + h M_L(\sigma)}\Bigr) & \le 0 \\
\Leftrightarrow & \frac{\p^3 }{\p h^3} \log \Bigl(  \frac {\sum e^{-\beta_c E_L + h M_L}} {\sum e^{-\beta_c E_L}} \Bigr) 
& \le 0 \\
\Leftrightarrow &  \frac{\p^2}{\p h^2} \Bigl(
\frac
{\frac {\p} {\p h} \<{e^{h M_L}}_{\beta_c,0,+}}
{\<{e^{h M_L}}_{\beta_c,0,+}}
 \Bigr) & \le 0\,.
\end{eqnarray*}
%where $E_L=E_L(\sigma)$ denotes the term in the Hamiltonian corresponding to $-\sum_{i\sim j} \sigma_i \sigma_j$. 
\ni
 
Let $F(h) = F_L(h) :=  \frac 
{\frac \p {\p h} \<{e^{h M_L}}_{\beta_c,0,+}}
{\<{e^{h M_L}}_{\beta_c,0,+}} = \<{M_L}_{\beta_c,h,+}$.
Then one has for any $h\geq 0$:
\begin{align}\label{e.1and2}
F(h) &\le F(0) + h\, F'(0) \nn\\
& = \<{M_L}_{\beta_c,0,+} + h \bigl( \<{M_L^2}_{\beta_c,0,+} - \<{M_L}_{\beta_c,0,+}^2 \bigr)\,.
% \\
%& \le C \bigl( L^{15/8} + h\, L^{15/4} \bigr) \,,
\end{align}
We will use the following Proposition from \cite{\CGNSL}, whose proof relies essentially on the RSW
theorem for Ising-FK percolation proved in \cite{\RSWfk}.

\begin{proposition}[Proposition B.2 in \cite{\CGNSL}] 
\label{pr.1and2}
There is a universal constant $C>0$ such that for $L$ sufficiently large, one has 
\bi
\item[(i)] $\<{M_L}_{\beta_c,0,+} \le C \, L^{2} \rho(L)^{1/2}$
\item[] and 
\item[(ii)] $\<{M_L^2}_{\beta_c,0,+}  \le C \, L^4 \, \rho(L)$\,.
\ei
\end{proposition}

\medskip
Plugging the inequalities of Proposition \ref{pr.1and2} into~\eqref{e.1and2} and using Wu's result, Theorem \ref{th.Wu}, gives us the following upper bound for $F(h)$:
\begin{equation}\label{}
F(h) \le C \bigl( L^{15/8} + h\, L^{15/4} \bigr) \,.
\end{equation}

Plugging this into~\eqref{e.dombymagn} gives us
\begin{equation} \label{upper-bound}
\<{\sigma_0}_{\beta_c, h} \le \frac 1 {L(h)^2}  \, F(h) \le C(L^{15/8} + h L^{15/4})/L^2 \,.
\end{equation}
Optimizing in $L = L(h) \geq 1$ leads to $ h\, L(h)^{15/4} \asymp L(h)^{15/8}$, which in turn gives
\[
L(h)\asymp h^{-8/15}\,.
\]
Thus
\begin{align}
\<{\sigma_0}_{\beta_c, h} & \le O(1) \, \frac  1 {L(h)^2} L(h)^{15/8}  = O(1)\, L(h)^{-1/8}\\
& \le O(1) \, h^{1/15}\,.
\end{align}
which concludes the proof of the upper bound in Theorem \ref{th.hZ2}.
\QED

\begin{remark}
We note that optimizing inequality \eqref{upper-bound}, which leads to $L(h)\asymp h^{-8/15}$, is
equivalent to letting $L(h)$ scale like the correlation length of the  near-critical model for vanishing external magnetic
field (i.e., like the inverse of the so-called ``mass scale'').  Indeed, the correlation length, $\xi$, scales like
$\xi(h) \asymp h^{-8/15}$.
\end{remark}

\begin{remark}
One might wonder why dominating by the $+$ boundary condition was a necessary step here.
The reason is that one cannot apply the GHS inequality to the finite volume quantity $M_L^{\Z^2}$, since the magnetization field 
is also increased ``outside'' $\Lambda_L$. Note that if one had worked without the $+$ boundary condition, then $F(0)$ would be zero, 
and one would obtain a suspicious upper bound on $\<{\sigma_0}_{\beta_c, h}$ of the form $h \Var{M_L}_{\beta_c,0} \asymp h\, L^{15/4}$ valid for all $L\geq 1$. This would clearly be wrong since, taking $L=1$, it would give an upper bound on $\<{\sigma_0}_{\beta_c, h}$ of order $h$, which is too small compared to the correct behavior of $h^{1/15}$.
\end{remark}

%\QED
\medskip

%\ni
%\underline{Lower bound.}
\section{Proof of the lower bound}\label{s.lower}

%The lower bound can be obtained combining the Buckingham-Gunton inequality \cite{GuBu68,Fisher69}
%with Wu's result, Theorem \ref{th.Wu}, and RSW. %(the latter is needed to obtain $\eta=1/4$ in all directions).
%%For completeness, 
%We provide here an elementary, self-contained proof based on RSW. \margin{F: new sentence}

%(We also remark that the approach of \cite{Newman79} to the BG inequality, combined with T. T. Wu's result
%and RSW, yields a lower bound for the free energy, $\int_0^h <\sigma_0>_{\beta,h} dh$, of the
%form $c h^{(1+1/15)}$.)

%The lower bound follows from rather standard percolation-type arguments, using the ghost spin representation
%of the Ising model with an external magnetic field and the RSW theorem of \cite{arXiv:0912.4253}.
%%The lower bound is more classical and may have been already known although we did not find a reference. 
%For completeness, we include it here. 
%%Here is an outline of how it works.
%We also remark that the lower bound is related to the Buckingham-Gunton
%inequality \cite{GuBu68,Fisher69}. In particular, the approach to that inequality
%of \cite{Newman79} can be combined with T. T. Wu's result and RSW to obtain a lower
%bound for the free energy, $\int_0^h <\sigma_0>_{\beta,h} dh$, of the
%form $c h^{(1+1/15)}$. This lower bound does not seem to imply the bound we are looking for. 
%\medskip

The lower bound is related to the Buckingham-Gunton inequality \cite{GuBu68,Fisher69}.
In particular, the approach to that inequality of \cite{Newman79} can be combined with T. T. Wu's result
and RSW to obtain a lower bound for the free energy, $\int_0^h <\sigma_0>_{\beta,h} dh$, of the
form $c h^{(1+1/15)}$. %Yet this lower bound does not seem to easily imply the bound we are looking for. 
%We provide here an elementary, self-contained proof based on the  RSW theorem from \cite{\RSWfk}.
Here we provide a lower bound for the magnetization, matching the upper bound of the previous section.
The proof is elementary and self-contained and is based on standard percolation arguments, using the
ghost spin representation of the Ising model with an external magnetic field and the RSW theorem of \cite{\RSWfk}.
\medskip

Let us first settle some notation.
We denote by $\P_{p_c,h}$ the $\fk$ percolation model representing the Ising model at $\beta_c$ with 
positive magnetic field $h\geq 0$. It is  a model of FK percolation defined on the extended graph $G=(V,E)$ where $V=\Z^2 \cup \{\g\}$ (the vertex $\g$ is commonly called the {\bf ghost vertex}) and $E$ is the usual set of nearest-neigbor edges $\E^2$ plus all the edges of the form $e= \{ x,\g \}$. %$e=\<{x,\g}$. 
(See, for example, \cite{\Griffith,\ABF87,\Chayes} for the use of the ghost vertex in the Ising model and in FK percolation.)
Furthermore, each edge $e\in \E^2$ carries a weight $p_c = 1-e^{-\beta_c}$, while each edge $e=\{x,\g\}$ carries a weight $p_h:= 1- e^{-h}$. 
It is well known that $\P_{p_c, h}$ stochastically dominates $\P_{p_c}$ (denoted $\P_{p_c, h} \succcurlyeq \P_{p_c}$). 
Furthermore, in the plane $\Z^2$, one has
\begin{align*}
\<{\sigma_0}_{\beta_c,h} & = \FK{p_c, h}{}{0 \leftrightarrow \g}\, ,
\end{align*}
where $0 \leftrightarrow \g$ denotes the event that the edge between the origin, $0$, and the ghost vertex, $\g$,
is open.

In the rest of the paper, the $\fk$ percolation configurations sampled according to $\P_{p_c,h}$ will be denoted 
$\bar \omega_h = (\omega_h, \tau_h)$,
where the component $\omega_h$ denotes the configuration of open edges lying in $\Z^2$, while the component $\tau_h$ corresponds to the set of open edges 
going from points in $\Z^2$ to $\g$.

For any configuration of edges $\omega$ in the plane (i.e., any $\omega \subset \mathbb{E}^2$), 
let $\cluster(\omega)$ be the connected component of the origin. 
Note that with such a definition, $\cluster(\omega_h)$ might be strictly smaller than the connected component of
the origin in the (enlarged) configuration $\bar \omega_h$.

Since $\P_{p_c,h}  \succcurlyeq  \P_{p_c} $, the restriction $\omega_h$ stochastically dominates the standard
$\omega_0 \sim \P_{p_c,0}$. 
%\margin{It is fine since $\P_{p_c,h}$ is FK percolation on a 
%larger graph than $\P_{p_c}$.} 
In particular, 
\[
\<{|\cluster(\omega_h)|}_{p_c,h} \geq \<{|\cluster(\omega_0)|}_{p_c,0} \,.
\]

For $M\in \N^+$, let $A_M$ be the event 
\begin{align*}
A_M(\omega):= \{ |\cluster(\omega)| \geq M\}\,.
\end{align*}

For any $M\in \N^+$, one has
\begin{align}\label{e.lowerbdsigma}
\<{\sigma_0}_{\beta_c,h} & = \FK{p_c, h}{}{0 \leftrightarrow \g}  \nn \\
& \geq  \FK{p_c, h}{}{A_M(\omega_h) \text { and } 0 \leftrightarrow \g } \nn  \\
& =  \FK{p_c, h}{}{A_M(\omega_h)} \, \FK{p_c, h}{}{0 \leftrightarrow \g \md A_M }  \nn  \\
&\geq  \FK{p_c, 0}{}{A_M(\omega_0)} \, \FK{p_c, h}{}{0 \leftrightarrow \g \md A_M }\,,
\end{align} 
since the event $A_M$ is clearly increasing and since $\omega_h \succcurlyeq  \omega_0$.  

In order to conclude the proof of the lower bound, we use the following two lemmas.

\begin{lemma}\label{l.M115}
There exists a constant $c_1>0$ such that for any $M\in \N^+$,
\begin{align*}
\FK{p_c,0}{}{A_M(\omega_0)} & = \FK{p_c, 0}{}{|\cluster(\omega_0)| \geq M} \geq c_1 \, M^{-1/15}\,.
\end{align*}

\end{lemma}

\begin{lemma}\label{l.upperboundghost}
There is a constant $c_2>0$ such that  for any $M\in \N^+$ and any $h \geq 0$,
\[
\FK{p_c, h}{}{0 \nleftrightarrow \g \md  A_M } \leq e^{-c_2 \,h \, M}\,.
\]
\end{lemma}

Before detailing the proofs of these lemmas, let us see why they enable us to conclude the proof
of the lower bound.
By combining the above lemmas with~\eqref{e.lowerbdsigma}, one has that, for any $M\in \N^+$,

\begin{align*}
\<{\sigma_0}_{\beta_c,h} & \geq \Omega(1) \, M^{-1/15}\, \bigl( 1 - e^{-c_2 \,h \, M} \bigr)\,.
\end{align*}

Now, optimizing $M\in \N^+$, one finds that $M$ should be chosen so that $M\asymp h^{-1}$. 
This particular choice of $M$ gives the lower bound of $\Omega(1) h^{1/15}$. \QED

It remains to prove Lemma \ref{l.M115} and Lemma \ref{l.upperboundghost}.

\medskip
\ni
{\bf Proof of Lemma \ref{l.M115}.}
%\medskip
\ni
%%%%%%%%%%%%%%%%%%%%%%%%%%%
%Let $R:= \lambda\, M^{8/15}$ ($\lambda$ is a fixed constant which will  be determined later) and let $C_R$ be the event that there is an open circuit in the annulus $A(R/2,R)$.
For any radius $R\geq 1$, let $C_R$ be the event that there is an open circuit in the annulus
$A(R/2,R) := \Lambda_R \setminus \Lambda_{R/2}$.
We will also write $0 \leftrightarrow R$ to denote the event,  $0 \leftrightarrow \partial\Lambda_R$,
that the origin is connected to $\partial \Lambda_R$ by a path of open edges.

\begin{align}
\FK{p_c,0}{}{A_M} & = \FK{p_c, 0}{}{|\cluster(\omega_0)| \geq M} \nonumber \\
& \geq  \FK{p_c, 0}{}{0 \leftrightarrow R,\, C_R } \FK{p_c, 0}{}{|\cluster(\omega_0)| \geq M \md 0 \leftrightarrow R,\, C_R} \nonumber \\
& \geq \FK{p_c, 0}{}{0 \leftrightarrow R }\, \FK{p_c, 0}{}{ C_R } \, \FK{p_c, 0}{}{|\cluster(\omega_0)| \geq M \md 0 \leftrightarrow R,\, C_R}  \;\; \text{(by FKG)} \nonumber \\
& \geq \Omega(1) R^{-1/8} \; \FK{p_c, 0}{}{|\cluster(\omega_0)| \geq M \md 0 \leftrightarrow R,\, C_R}\,, \label{e.MM}
%& \geq \Omega(1)  M^{-1/15} 
\end{align}
by RSW from \cite{\RSWfk} (see also Lemma  B.3 of \cite{CGNSL}).

It remains to prove that if the % constant $\lambda>0$ is chosen large enough, then 
radius $R\geq 1$ is chosen (as a function of $M$) to be $\Omega(1) M^{8/15}$, then one has 

\[
\FK{p_c, 0}{}{|\cluster(\omega_0)| \geq M \md 0 \leftrightarrow R,\, C_R} \geq \Omega(1)\,.
\]

This is easily done by a second moment argument on the random variable $N:= |\cluster(\omega_0) \cap \Lambda_R| $.
%More precisely, 
Indeed, denoting $x+\Lambda_R$ by $B(x,R)$, we have
\begin{align*}
\EFK{p_c,0}{}{N  \md 0 \leftrightarrow R,\, C_R} & = \sum_{x\in \Lambda_R} \FK{p_c, 0}{}{  0 \leftrightarrow x  \md 0 \leftrightarrow R,\, C_R} \nn \\
& \geq \sum_{x\in B(0,R/2)} \FK{p_c, 0}{}{   x \leftrightarrow \p B(x,2R) \md 0 \leftrightarrow R,\, C_R} \nn \\
& \geq  \Omega(1) R^2 \FK{p_c, 0}{}{  0 \leftrightarrow 2R}  \;\;  \left( \begin{array}{l} \text{By translation inva-} \\ \text{riance and FKG} \end{array} \right) \nn \\
& \geq \Omega(1) R^{15/8}\,. %\label{e.MM}
\end{align*}

For the second moment, one has 
\begin{align*}
\EFK{p_c,0}{}{N^2 \md 0 \leftrightarrow R,\, C_R} &=  \sum_{x,y \in \Lambda_R} \FK{p_c, 0}{}{  0 \leftrightarrow x, 0\leftrightarrow y  \md 0 \leftrightarrow R,\, C_R} \\
& \le  \sum_{x,y \in \Lambda_R} \frac  {\FK{p_c, 0}{}{  0 \leftrightarrow x, 0\leftrightarrow y ,\, 0 \leftrightarrow R,\, C_R}} {\FK{p_c,0}{}{0 \leftrightarrow R}\, \FK{p_c,0}{}{C_R} }  \; \text{(by FKG)} \\
& \le O(1) R^{1/8}  \sum_{x,y \in \Lambda_R}  \,  \FK{p_c,0}{}{ 0 \leftrightarrow x, 0\leftrightarrow y ,\, 0 \leftrightarrow R} \,.
\end{align*}

Now, a standard dyadic summation (as for example in Section 4.3 or 7.2 of \cite{\GPS}) gives 

\[
\sum_{x,y \in \Lambda_R} \FK{p_c,0}{}{ 0 \leftrightarrow x, 0\leftrightarrow y ,\, 0 \leftrightarrow R}  \asymp R^{4} \, \FK{p_c, 0}{}{0 \leftrightarrow R}^3 \, ,
\]

which implies the following upper bound on the second moment:
\[
\EFK{p_c,0}{}{N^2 \md 0 \leftrightarrow R,\, C_R} \le O(1) R^{15/4}\,.
\]

By the second moment method, there exists a constant $c>0$ such that 
\begin{equation}\label{e.MM2}
\FK{p_c, 0}{}{|\cluster(\omega_0)| \geq c R^{15/8} \md  0 \leftrightarrow R,\, C_R  } > c \, .
\end{equation}

Now let us choose $R=R(M)$ so that $c R^{15/8} = M$, i.e., $R:= C M^{8/15}$.
Plugging this choice into ~\eqref{e.MM} and using ~\eqref{e.MM2} concludes the proof of Lemma \ref{l.M115}.

\QED
%%%%%%%%%%%%%%%%%%%%%%%%%%%

\ni
{\bf Proof of Lemma \ref{l.upperboundghost}.}
%\medskip
\ni
%%%%%%%%%%%%%%%%%%%%%%%%%%%
Sample the configuration $\omega_h$ (which by definition is the configuration $\bar \omega_h$ without the edges going from $\Z^2$ to the ghost $\bf g$) according to the conditional measure $\FK{p_c,h}{}{\cdot \md A_M}$. Let $k:= |\cluster(\omega_h)|$ (in particular $k\geq M$) and index the $k$ vertices 
in $\cluster(\omega_h)$ in any order: $x_1, \ldots x_k$. Let $\tau_i$ for $i=1,\ldots,k$ be the edge going from vertex $x_i$ to the ghost.
Sample the edges $\tau_i$ one at a time knowing the configuration $\omega_h$ as well as  the edges  $\tau_j$ already sampled.
For any $i\in \{0,\ldots, k-1\}$, by a quantitative version of the {\bf finite energy property}, one has 
\[
\FK{p_c, h}{}{\tau_{i+1} \text{ is open} \md \omega_h,\, \tau_1,\ldots, \tau_i} \geq c\, h\,,
\]
for some constant $c>0$. 
This implies that 
\begin{align*}
\FK{p_c, h}{}{0  \nleftrightarrow {\bf g}  \md A_M} & = \sum_{\omega_h} \FK{p_c,h}{}{\omega_h \md A_M}\, \FK{p_c, h}{}{\cluster(\omega_h) \nleftrightarrow {\bf g}  \md \omega_h,\, A_M}  \\
& \le (1-ch)^{|\cluster(\omega_h)|} \\
&\le (1- ch)^M\,,
\end{align*}
which easily concludes the proof of Lemma \ref{l.upperboundghost}.
\QED

\begin{remark}\label{}
One may wonder which parts of the proof (upper or lower bound) use the planarity of the model. First of all,
the main ingredient used for the upper bound, i.e.,
the GHS inequality, is of course valid in any dimension. Yet,
planarity is used for the upper bound since the proof uses the fact that the variance of the magnetization
behaves like the first moment squared. This fact is the content of Proposition B.2 in \cite{\CGNSL}, which
relies on the RSW Theorem of \cite{\RSWfk} and thus relies in an essential way on planarity. Our proof of
the lower bound also relies on the RSW theorem of \cite{\RSWfk} and thus also requires planarity.
\end{remark}

\section{Without assuming T. T. Wu's result}\label{WWu}
If one does not want to assume Wu's result, instead of Theorem \ref{th.hZ2}
one obtains the following result.
\begin{theorem}[without assuming T. T. Wu's result]\label{}
Consider the Ising model on $\Z^2$ at $\beta_c$ with a positive external magnetic field 
$h>0$, then 
\[
\<{\sigma_0}_{\beta_c, h} \asymp \sqrt{ \rho \left( \xi(h) \right) } \asymp \frac{1}{h\, \xi(h)^2}\,,
\]
where the {\bf correlation length} $\xi(h)$ is defined as follows:
\begin{equation}\label{e.CL}
\xi(h):= \inf  \left\{ L\geq 1:\, L^2 \sqrt{ \rho \left( L \right) } \geq \frac 1 h \right\} \,.
\end{equation}
\end{theorem}

\begin{remark}\label{r.SLE}
Note that if one could show that %$\alpha_a^\fk(\eps, 1) \asymp \eps^{1/8}$
$\alpha_a^\fk(\eps, 1) := \FK{p_c,0}{}{\partial\Lambda(\epsilon/a) \leftrightarrow \partial\Lambda(1/a)} \asymp \eps^{1/8}$
with an $SLE$ computation,
this would imply, without assuming T. T. Wu's result, the following estimate on the average magnetization at
$\beta= \beta_c$ as $h \to 0+$:
\[
\<{\sigma_0}_{\beta_c, h} = h^{1/15+o(1)}\,.
\]
This highlights that we do not need the full strength of T. T. Wu's result to obtain the exponent~$\frac 1 {15}$. 
\end{remark}

\medskip
\ni
{\bf Proof.}
The proof of the upper bound works exactly in the same fashion as before since one can rely on Proposition \ref{pr.1and2} which was proved in the appendix of \cite{\CGNSL} without relying on Wu's result. 

For the proof of the lower bound, we need to replace Lemma \ref{l.M115} by the following (note that the proof of Lemma \ref{l.upperboundghost} did not assume Wu's result).
 
\begin{lemma}\label{l.M115B}
There exists a constant $c_1>0$ such that for any $M\in \N^+$,
\begin{align*}
\FK{p_c,0}{}{A_M} & = \FK{p_c, 0}{}{|\cluster(\omega_0)| \geq M} \geq c_1 \, \sqrt{\rho\big( \xi( c_1 / M)\big)}\,,
\end{align*}
where $\xi(\cdot)$ is the correlation length defined in~\eqref{e.CL}.
\end{lemma}

To see why this holds, one proceeds in the same fashion as in the proof of Lemma \ref{l.M115}, using some radius $R$ to be chosen later (as a function of $M$). 
As in the appendix of \cite{\CGNSL}, one finds the following bounds on the first and second moments of the random variable $N$ defined above:
\begin{equation*}
\begin{cases}
\EFK{p_c,0}{}{N  \md 0 \leftrightarrow R,\, C_R} \geq \Omega(1) \, R^2 \rho(R)^{1/2} \\
\EFK{p_c,0}{}{N^2  \md 0 \leftrightarrow R,\, C_R}  \le O(1) R^4 \rho(R)
\end{cases}
\end{equation*}

Therefore by the second moment method, there is some $c>0$ such that $N$ is larger than $c R^2 \rho(R)^{1/2}$ with positive conditional probability at least $c>0$.
Now with 
\begin{equation}\label{}
R=R(M):= \inf \{ R \geq 1, R^2 \rho(R)^{1/2} \geq \frac 1 c M\}=: \xi(c/M)\,,
\end{equation}
the same proof as above gives us that 
\begin{align*}\label{}
\FK{p_c,0}{}{A_M}  &\geq c\,  \FK{p_c, 0}{}{0 \leftrightarrow \xi(c/M) }  \\
& \geq \frac c {C} \sqrt{ \rho\big( \xi(c/M)\big)} \,,
\end{align*}
where we use in the last inequality the bound on $\FK{p_c, 0}{\free}{0 \leftrightarrow N} $ from Lemma A.3 of \cite{\CGNSL}.
This concludes our proof with $c_1 = c/C<c$. \QED

Combining the above estimates, and using $M:= \frac {c_1} h$,  we find that 
\begin{align}%\label{e.lowerbdsigma}
\<{\sigma_0}_{\beta_c,h} & = \FK{p_c, h}{}{0 \leftrightarrow {\bf g}}  \nn \\
&\geq  \FK{p_c, 0}{}{A_M(\omega_0)} \, \FK{p_c, h}{}{0 \leftrightarrow {\bf g} \md A_M } \\
&\geq  c_1 \sqrt{ \rho\big( \xi(h) \big)} (1- ch)^{c_1/h}  \\
&\geq \Omega(1) \sqrt{ \rho\big( \xi(h) \big)}\,,
\end{align} 
which completes the proof of Theorem \ref{WWu}. \QED

\medskip
\ni
{\bf Acknowledgments.} The authors thank Douglas Abraham, Hugo Duminil-Copin and
Roberto Fern\'andez for useful discussions,  and an anonymous referee for useful comments.

\ \\
{\bf Federico Camia}\\
Department of Mathematics, Vrije Universiteit Amsterdam and NYU Abu Dhabi \\
Research supported in part by NWO grant Vidi 639.032.916.\\
\\
{\bf Christophe Garban}\\
ENS Lyon, CNRS\\
%\url{http://www.umpa.ens-lyon.fr/~cgarban/}\\
Partially supported by ANR grant BLAN06-3-134462.\\
\\
{\bf Charles M. Newman}\\
Courant institute of Mathematical Sciences, New York University, New York, NY 10012, USA\\
Research supported in part by NSF grants OISE-0730136 and MPS-1007524.

\end{document}